\documentclass[10pt]{article}

\usepackage{amssymb}
\usepackage{eucal}
\usepackage{amsmath}
\usepackage{amsthm}
\usepackage{amscd}
\usepackage{graphics}
\usepackage{graphicx}
\usepackage{amsfonts}
\usepackage{latexsym}
\usepackage[all]{xy}

\numberwithin{equation}{section}
\newtheorem{thm}{\bf Theorem}[section]
\newtheorem{lem}[thm]{\bf Lemma}

\theoremstyle{remark}

\begin{document}

\title{Energy methods in the stability problem for the $\mathfrak{so}(4)$ free rigid body}
\author{Petre Birtea and Ioan Ca\c{s}u\\
{\small Department of Mathematics, West University of Timi\c soara,}\\
{\small Bd. V. Parvan, No. 4, 300223 Timi\c soara, Romania}\\
{\small E-mail: birtea@math.uvt.ro; casu@math.uvt.ro}}
\date{}

\maketitle

\begin{abstract}
For the $\mathfrak{so}(4)$
free rigid body the stability problem for isolated equilibria has been completely solved using Lie-theoretical and topological arguments. For each case of nonlinear stability previously found we construct a Lyapunov function. These Lyapunov functions are linear combinations of Mishchenko's constants of motion.  

\end{abstract}



\section{Introduction}

The goal of this paper is to prove stability using energy methods for the Cartan type equilibria of $\mathfrak{so}(4)$
free rigid body. More precisely, for each case of nonlinear stability, we will explicitly construct a Lyapunov function. We find these Lyapunov functions to be linear combinations of four integrals of motion which make the system completely integrable. These linear combinations of constants of motion are distinct for different nonlinear stable equilibria.   

As proved in \cite{noi}, there are three Cartan families of equilibria, organized in Weyl group orbits ge\-ne\-ra\-ted by three coordinate type Cartan subalgebras intersected with a regular adjoint orbit. Their
stability is studied in \cite{noi} using exclusively a Lie-theoretical result due to Williamson \cite{will}. The above
equilibria are grouped in three categories:
\begin{itemize}
\item[(I)] a class of center-saddle unstable equilibria;
\item[(II)] a class of center-center stable equilibria for which the stability does not depend on the parameters
defining the corresponding orbit;
\item[(III)] a class of equilibria for which a bifurcation phenomenon occurs, depending on the parameters
defining the corresponding orbit (the stability region being also of center-center type).
\end{itemize}

For equilibria of type (II) we find Lyapunov functions which guarantee their stability using linear combinations of Mishchenko's constants of motion, as described in \cite{Popov}. 

For equilibria of type (III) we first study their spectral stability and we find that a bifurcation phenomenon in the spectral stability occurs. The region of spectral stability overlaps the region of nonlinear stability previously discovered in \cite{noi}. For each region of nonlinear stability found in \cite{noi} we construct a Lyapunov function. We also point out a frontier case of spectral stability which is nonlinear unstable. 

Although the constants of motion used in our paper are functionally dependent of the Mishchenko's first integrals,
they prove their utility in applying energy methods for studying stability of equilibria. This paper completes previous results on nonlinear stability using energy methods \cite{FeMa03}, \cite{Spiegler04}. 
The stability problem and the implication of the topological structure of the energy-momentum level sets on bifurcation phenomena in the dynamics of the $\mathfrak{so}(4)$ free rigid body was extensively studied in \cite{Oshemkov87}, \cite{Oshemkov91}, \cite{Pogosyan83}, \cite{Pogosyan84}, \cite{tonkonog}.

We prove that Lyapunov functions are positive or negative definite studying the non-degeneracy of the Hessian matrix associated to the Lyapunov function at the equilibrium point. If the Hessian is degenerate at an equilibrium point, one can still study the Lyapunov stability using algebraic methods, see \cite{dan}, \cite{dan2}, \cite{dan3}.

\section{The geometry and stability problem for the $\mathfrak{so}(4)$ free rigid body}

In the first part of this section, we recall the geometry underlying the system of the free rigid body on the Lie algebra $\mathfrak{so}(4)$, as it was presented in \cite{noi}.
The equations of the rigid body on $\mathfrak{so}(4)$ are given by
\begin{equation}\label{rigideq}
\dot M=[M,\Omega], \end{equation} where $\Omega\in \mathfrak{so}(4)$,
$M=\Omega J+J\Omega\in \mathfrak{so}(4)$ with $J=\operatorname{diag}(\lambda_i)$, a real constant diagonal matrix satisfying $ \lambda _i + \lambda _j \geq 0 $, for all $i, j=1, \ldots, 4$, $i \neq j $ (see, for example, \cite{Ratiu80}).
Note that $M=[m_{ij}]$ and $\Omega=[ \omega_{ij}]$ determine each other if and only if $\lambda_i + \lambda _j >0$ since $m_{ij} = ( \lambda_i+ \lambda_j) \omega_{ij}$ which physically means that the rigid body is not concentrated on a lower dimensional subspace of $\mathbb{R}^6$. 

It is well known that equations \eqref{rigideq} are Hamiltonian relative to the minus Lie-Poisson bracket on $\mathfrak{so}(4)$
\begin{equation}
\label{LP_n}
\{F,G\}(M) : = \frac{1}{2}\operatorname{Trace}(M [\nabla F(M), \nabla G(M)]),
\end{equation}
and the Hamiltonian function
\begin{equation}
\label{ham_n}
H(M) := - \frac{1}{4} \operatorname{Trace}(M \Omega).
\end{equation}
We choose as basis of $\mathfrak{so}(4)$ the matrices
$$E_1=\left[\begin{array}{cccc}
0&0&0&0\\
0&0&-1&0\\
0&1&0&0\\
0&0&0&0\end{array}\right];~E_2=\left[\begin{array}{cccc}
0&0&1&0\\
0&0&0&0\\
-1&0&0&0\\
0&0&0&0\end{array}\right];~E_3=\left[\begin{array}{cccc}
0&-1&0&0\\
1&0&0&0\\
0&0&0&0\\
0&0&0&0\end{array}\right];$$
$$E_4=\left[\begin{array}{cccc}
0&0&0&1\\
0&0&0&0\\
0&0&0&0\\
-1&0&0&0\end{array}\right];~E_5=\left[\begin{array}{cccc}
0&0&0&0\\
0&0&0&1\\
0&0&0&0\\
0&-1&0&0\end{array}\right];~E_6=\left[\begin{array}{cccc}
0&0&0&0\\
0&0&0&0\\
0&0&0&1\\
0&0&-1&0\end{array}\right]
$$
and hence we represent $\mathfrak{so}(4)$ as
\begin{equation}
\label{so4_representation}
\mathfrak{so}(4)=\left\{\left.
M=\left[\begin{array}{cccc}
0&-x_3&x_2&y_1\\
x_3&0&-x_1&y_2\\
-x_2&x_1&0&y_3\\
-y_1&-y_2&-y_3&0\end{array}\right] \, \right|\, x_1,x_2,x_3,y_1,y_2,y_3\in \mathbb{R}
\right\}.
\end{equation}
Since $\hbox{rank}\,\mathfrak{so}(4) = 2$, we have two functionally independent Casimir functions which are given respectively by
\[
C_1(M) := -\frac{1}{4}\operatorname{Trace}(M^2)
= \frac{1}{2}\left(\sum\limits\limits_{i=1}^3 x_i^2+\sum\limits\limits_{i=1}^3 y_i^2 \right)
\]
and
\[
C_2(M) := - \operatorname{Pf}(M)
=\sum\limits\limits_{i=1}^3 x_iy_i.
\]
Thus the generic adjoint orbits are the level sets
$$
\operatorname{Orb}_{c_1c_2}(M)=(C_1\times C_2)^{-1}(c_1,c_2), 
\qquad (c_1,c_2) \in \mathbb{R}^2.
$$

In all that follows we will denote by $\operatorname{Orb}_{c_1;c_2}$ the regular adjoint orbit $\operatorname{Orb}_{c_1c_2}$,
where $c_1>0$ and $c_1>|c_2|$.

We will work from now on with a generic $\mathfrak{so}(4)$-rigid body, that is, $\lambda_i+ \lambda_j>0$ for $i \neq j$ and all $\lambda_i $ are distinct. The relationship between $\Omega=[\omega_{ij}]\in \mathfrak{so}(4)$ and the matrix $M\in \mathfrak{so}(4)$ in the representation \eqref{so4_representation} is hence given by
$$\left.\begin{array}{lll}
(\lambda_3+\lambda_2)\omega_{32}=x_1& \quad 
(\lambda_1+\lambda_3)\omega_{13}=x_2& \quad
(\lambda_2+\lambda_1)\omega_{21}=x_3\\
(\lambda_1+\lambda_4)\omega_{14}=y_1& \quad 
(\lambda_2+\lambda_4)\omega_{24}=y_2& \quad 
(\lambda_3+\lambda_4)\omega_{34}=y_3
\end{array}\right.$$
and thus the equations of motion \eqref{rigideq} are equivalent to the system
\begin{equation}\label{soM}
\left\{\begin{array}{l}
\vspace{.1cm}
\dot x_1=\left(\frac{1}{\lambda_1+\lambda_2}-\frac{1}{\lambda_1+\lambda_3}\right)x_2x_3+\left(\frac{1}{\lambda_3+\lambda_4}-\frac{1}{\lambda_2+\lambda_4}\right)y_2y_3\\
\vspace{.1cm}
\dot x_2=\left(\frac{1}{\lambda_2+\lambda_3}-\frac{1}{\lambda_1+\lambda_2}\right)x_1x_3+\left(\frac{1}{\lambda_1+\lambda_4}-\frac{1}{\lambda_3+\lambda_4}\right)y_1y_3\\
\vspace{.1cm}
\dot x_3=\left(\frac{1}{\lambda_1+\lambda_3}-\frac{1}{\lambda_2+\lambda_3}\right)x_1x_2+\left(\frac{1}{\lambda_2+\lambda_4}-\frac{1}{\lambda_1+\lambda_4}\right)y_1y_2\\
\vspace{.1cm}
\dot y_1=\left(\frac{1}{\lambda_3+\lambda_4}-\frac{1}{\lambda_1+\lambda_3}\right)x_2y_3+\left(\frac{1}{\lambda_1+\lambda_2}-\frac{1}{\lambda_2+\lambda_4}\right)x_3y_2\\
\vspace{.1cm}
\dot y_2=\left(\frac{1}{\lambda_2+\lambda_3}-\frac{1}{\lambda_3+\lambda_4}\right)x_1y_3+\left(\frac{1}{\lambda_1+\lambda_4}-\frac{1}{\lambda_1+\lambda_2}\right)x_3y_1\\
\vspace{.1cm}
\dot y_3=\left(\frac{1}{\lambda_2+\lambda_4}-\frac{1}{\lambda_2+\lambda_3}\right)x_1y_2+\left(\frac{1}{\lambda_1+\lambda_3}-\frac{1}{\lambda_1+\lambda_4}\right)x_2y_1.
\end{array}\right.
\end{equation}
The Hamiltonian \eqref{ham_n} has the expression
\begin{align*}
H(M)&=-\frac{1}{4}\hbox{Trace}(M\Omega)\\
&=\frac{1}{2}\left(\frac{1}{\lambda_2+\lambda_3}x_1^2+\frac{1}{\lambda_1+\lambda_3}x_2^2+\frac{1}{\lambda_1+\lambda_2}x_3^2+\frac{1}{\lambda_1+\lambda_4}y_1^2+\frac{1}{\lambda_2+\lambda_4}y_2^2+\frac{1}{\lambda_3+\lambda_4}y_3^2\right).
\end{align*}

\begin{thm} \label{echilibre} 
\cite{noi} If $\mathcal{E}$ denotes the set of the equilibrium points
of \eqref{soM}, then
${\cal E}=\mathfrak{t}_1\cup \mathfrak{t}_2\cup \mathfrak{t}_3\cup \mathfrak{s}_+ \cup \mathfrak{s}_-$,
where
$$\mathfrak{t}_1:
=\left\{\left. M_{a,b}^1:=\left[\begin{array}{cccc}
0&0&0&b\\
0&0&-a&0\\
0&a&0&0\\
-b&0&0&0\end{array}\right]\, \right| \,a,b\in \mathbb{R}
 \right\},$$
$$\mathfrak{t}_2:
=\left\{\left.  M_{a,b}^2:=\left[\begin{array}{cccc}
0&0&a&0\\
0&0&0&b\\
-a&0&0&0\\
0&-b&0&0\end{array}\right] \, \right| \, a,b\in \mathbb{R}
 \right\},$$
$$\mathfrak{t}_3:
=\left\{\left. M_{a,b}^3:=\left[\begin{array}{cccc}
0&-a&0&0\\
a&0&0&0\\
0&0&0&b\\
0&0&-b&0\end{array}\right] \, \right| \,a,b\in \mathbb{R}
 \right\}$$
are the coordinate type Cartan subalgebras of $\mathfrak{so}(4)$ and
$\mathfrak{s}_{\pm}$ are the three dimensional vector subspaces given by\\
\begin{center}
$\mbox{\fontsize{12}{12}\selectfont $
\mathfrak{s}_{\pm}:=\operatorname{span}_{\mathbb{R}}\left\{\left(\frac{1}{\lambda_1+\lambda_4}E_1 \pm \frac{1}{\lambda_2+\lambda_3}E_4\right),\,\left(\frac{1}{\lambda_2+\lambda_4}E_2\pm \frac{1}{\lambda_1+\lambda_3}E_5\right), \,\left(\frac{1}{\lambda_3+\lambda_4}E_3\pm \frac{1}{\lambda_1+\lambda_2}E_6\right) \right\}.$}$
\end{center}
\end{thm}

The intersection of a regular adjoint orbit and a coordinate
type Cartan subalgebra has four elements which represents a Weyl group orbit. Thus we expect twelve isolated equilibria for the rigid body equations \eqref{rigideq}. Specifically, we have the following result.

\begin{thm} \label{Weil_group_orbit}
\cite{noi} The following equalities hold:
\begin{itemize}
\item[{\rm (i)}] $\mathfrak{t}_1\cap \operatorname{Orb}_{c_1;c_2}=\left\{
M_{a,b}^1,M_{-a,-b}^1,M_{b,a}^1,M_{-b,-a}^1
\right\}$,
\item[{\rm (ii)}]  $\mathfrak{t}_2\cap \operatorname{Orb}_{c_1;c_2}=\left\{
M_{a,b}^2,M_{-a,-b}^2,M_{b,a}^2,M_{-b,-a}^2
\right\}$,
\item[{\rm (iii)}]  $\mathfrak{t}_3\cap
    \operatorname{Orb}_{c_1;c_2}=\left\{
    M_{a,b}^3,M_{-a,-b}^3,M_{b,a}^3,M_{-b,-a}^3
    \right\}$,
\end{itemize}
where
\begin{equation}
\label{values_a_b}
\left\{
\begin{aligned}
a&=\frac{1}{\sqrt{2}}\left(\sqrt{c_1+c_2}+\sqrt{c_1-c_2}\right) \\
b&=\frac{1}{\sqrt{2}}\left(\sqrt{c_1+c_2}-\sqrt{c_1-c_2}\right).
\end{aligned}
\right.
\end{equation}
\end{thm}

\medskip

Further we remind the main results obtained in \cite{noi} on the nonlinear stability of the
equilibrium states ${\cal E}\cap \operatorname{Orb}_{c_1;c_2}$ for the dynamics \eqref{soM} on a generic adjoint orbit.

Since the system \eqref{soM} on a generic adjoint orbit is completely integrable (\cite{BoFo04}, \cite{Fo88}, \cite{Mishchenko70}, \cite{MiFo78}), for the 
$\mathfrak{so}(4)$ free rigid body we have a supplementary constant of motion. Using Mishchenko's method \cite{Mishchenko70}, \cite{Ratiu80}, we obtain the following additional constant of the motion for the equations \eqref{soM} commuting with $H$:
$$
I(M)=(\lambda_2^2+\lambda_3^2)x_1^2+(\lambda_1^2+\lambda_3^2)x_2^2+(\lambda_1^2+\lambda_2^2)x_3^2+(\lambda_1^2+\lambda_4^2)y_1^2+(\lambda_2^2+\lambda_4^2)y_2^2+(\lambda_3^2+\lambda_4^2)y_3^2.
$$
Without loss of generality, we can choose an ordering for $\lambda_i$'s, namely
$$\lambda_1>\lambda_2>\lambda_3>\lambda_4.$$

The following results on nonlinear stability hold \cite{noi}.

\begin{thm}
\label{t1}
\cite{noi} The equilibria $M_{a,b}^1, M_{-a,-b}^1 \in 
\mathfrak{t}_1 \cap \operatorname{Orb}_{c_1; c_2}$ are non-degenerate of 
center-center type and therefore nonlinearly stable on the corresponding adjoint orbit and also nonlinearly stable
for the Lie-Poisson dynamics on $\mathfrak{so}(4)$.
\end{thm}

We denote by $\alpha_1,\alpha_2$ ($\alpha_1<\alpha_2$) the distinct real roots of the quadratic function
$$\widetilde{f}(t)=\widetilde{S}t^2+\widetilde{T}t+\widetilde{U},$$
where
$$\widetilde{S}=(\lambda_1^2-\lambda_4^2)^2>0;~~~\widetilde{U}=(\lambda_2^2-\lambda_3^2)^2>0;$$
$$\widetilde{T}=-2[(\lambda_1^2-\lambda_2^2)(\lambda_3^2-\lambda_4^2)+(\lambda_1^2-\lambda_3^2)(\lambda_2^2-\lambda_4^2)]<0.$$

\begin{thm}
\label{bifurcation}
\cite{noi} Under the hypothesis $\lambda_1^2+ \lambda_4^2\not= \lambda_2^2+\lambda_3^2$ the following hold:
\begin{itemize}
\item[{\rm (i)}] If $\frac{b^2}{a^2}\in [0,\alpha_1)$, then the equilibria $M_{b,a}^1, M_{-b,-a}^1 \in 
\mathfrak{t}_1 \cap \operatorname{Orb}_{c_1; c_2}$ are non-degenerate unstable of saddle-saddle type on the adjoint orbit determined by $a $ and $b $ and hence are also unstable for the Lie-Poisson dynamics on $\mathfrak{so}(4)$.
\item[{\rm (ii)}] If $\frac{b^2}{a^2}\in (\alpha_1,\alpha_2)$, then the equilibria $M_{b,a}^1, M_{-b,-a}^1 \in 
\mathfrak{t}_1 \cap \operatorname{Orb}_{c_1; c_2}$ are non-degenerate unstable of focus-focus type on the adjoint orbit determined by $a $ and $b $ and hence are also unstable for the Lie-Poisson dynamics on $\mathfrak{so}(4)$.
\item[{\rm (iii)}] If $\frac{b^2}{a^2}\in (\alpha_2,1)$, then the equilibria $M_{b,a}^1, M_{-b,-a}^1 \in 
\mathfrak{t}_1 \cap \operatorname{Orb}_{c_1; c_2}$ are non-degenerate stable of center-center type on the adjoint orbit determined by $a $ and $b $ and hence are also nonlinearly stable for the Lie-Poisson dynamics on $\mathfrak{so}(4)$.
\item[{\rm (iv)}] If $\frac{b^2}{a^2}=\alpha_1$, then the equilibria $M_{b,a}^1, M_{-b,-a}^1 \in 
\mathfrak{t}_1 \cap \operatorname{Orb}_{c_1; c_2}$ are degenerate and unstable on the adjoint orbit determined by $a $ and $b $ and hence are also unstable for the Lie-Poisson dynamics on $\mathfrak{so}(4)$.
\item[{\rm (v)}] If $\frac{b^2}{a^2}=\alpha_2$, then the equilibria $M_{b,a}^1, M_{-b,-a}^1 \in 
\mathfrak{t}_1 \cap \operatorname{Orb}_{c_1; c_2}$ are degenerate and the stability problem on the adjoint orbit determined by $a $ and $b $ remains open.
\end{itemize}
\end{thm}

\begin{thm}
\label{bifurcation2}
\cite{noi} Under the hypothesis $\lambda_1^2+ \lambda_4^2= \lambda_2^2+\lambda_3^2$ the following holds:
\begin{itemize}
\item[{\rm (i)}] If $\frac{b^2}{a^2}\in [0,\alpha_1)$, then the equilibria $M_{b,a}^1, M_{-b,-a}^1 \in 
\mathfrak{t}_1 \cap \operatorname{Orb}_{c_1; c_2}$ are non-degenerate unstable of saddle-saddle type.
\item[{\rm (ii)}] If $\frac{b^2}{a^2}\in (\alpha_1,1)$, then the equilibria $M_{b,a}^1, M_{-b,-a}^1 \in 
\mathfrak{t}_1 \cap \operatorname{Orb}_{c_1; c_2}$ are non-degenerate unstable of focus-focus type.  
\item[{\rm (iii)}] If $\frac{b^2}{a^2}=\alpha_1$, then the equilibria $M_{b,a}^1, M_{-b,-a}^1 \in 
\mathfrak{t}_1 \cap \operatorname{Orb}_{c_1; c_2}$ are degenerate and unstable.
\end{itemize}
Thus these equilibria are also unstable for the Lie-Poisson dynamics on $\mathfrak{so}(4)$.
\end{thm}

\begin{thm}\label{t3}
\cite{noi} All four equilibria in $\mathfrak{t}_3 \cap \operatorname{Orb}_{c_1; c_2}$ are non-degenerate of center-center type and therefore nonlinearly stable on the corresponding adjoint orbit. These equilibria are also nonlinearly stable
for the Lie-Poisson dynamics on $\mathfrak{so}(4)$.
\end{thm}

\begin{thm}\label{t2}
\cite{noi} All four equilibria in $\mathfrak{t}_2 \cap \operatorname{Orb}_{c_1; c_2}$ are non-degenerate, of center-saddle type and therefore unstable on the corresponding adjoint orbit. Thus these equilibria are also unstable for the Lie-Poisson dynamics on $\mathfrak{so}(4)$.
\end{thm}

\section{Nonlinear stability using energy methods for equilibria without bifurcation be\-ha\-vi\-or}

We study nonlinear stability using Arnold's method, which is equivalent \cite{BiPu07} with Energy-Casimir method \cite{holm} and with Ortega-Ratiu method \cite{ortega}, for the equilibria listed in Theorems \ref{t1} and \ref{t3}. The nonlinear stability of these equilibria has been proved in \cite{noi} using a Lie-theoretical result of Williamson. In this section we will find a Lyapunov function for these equilibria. 

For the equilibria $M_{a,b}^1$ and $M_{-a,-b}^1$ in $\mathfrak{t}_1 \cap \operatorname{Orb}_{c_1; c_2}$ the computations are more precisely as follows. 
Consider the smooth function $F_{mn}\in
C^{\infty}(\mathfrak{so}(4),\mathbb{R})$, where $m,n$ are real numbers
$$F_{mn}(M)=I(M)+m C_1(M)+n C_2(M).$$
Choosing $m,n$ such that $dF_{mn}(M_{a,b}^1)=0$ and taking into account that
$d^2F_{mn}(M_{a,b}^1)|_{W\times W}$ is indefinite, where
$W:=\ker dC_1(M_{a,b}^1)\cap \ker dC_2(M_{a,b}^1)$, Arnold's method shows its limitations when used with the constant of motion $I$. 

Using Arnold's method with the Hamiltonian of the system, i.e. an energy function of the form $F_{mn}(M)=H(M)+m C_1(M)+n C_2(M)$, we obtain that equilibria $M_{a,b}^1$ and $M_{-a,-b}^1$ are nonlinear stable under the sufficient conditions on $\lambda_i$ and respectively on $a,b$:
\begin{equation*}
\lambda_1+\lambda_4<\lambda_2+\lambda_3
\end{equation*}
and
\begin{equation*}
\frac{b^2}{a^2}>\frac{(\lambda_1+\lambda_4)^2}{(\lambda_2+\lambda_3)^2}.
\end{equation*}
These stability conditions were also obtained in \cite{daisuke}.
This stability result is weaker than Theorem \ref{t1}.

From the above considerations it becomes necessary to look for a Lyapunov function that involves all four constants of motion,
\begin{equation}
\label{energy}
F_{mn}^{\mu_1\mu_2}:=\mu_1 H+\mu_2 I+mC_1+nC_2,
\end{equation}
where $\mu_1,\mu_2,m,n\in \mathbb{R}$.

In \cite{Popov} it has been presented another set of constants of motion for the $\mathfrak{so}(4)$ rigid body. Following the same ideas as in \cite{casu}, we approach the stability problem by energy methods using these constants of motion. They have simple and elegant expressions:

\medskip
$$G_1(M)=\frac{x_2^2}{\lambda_1^2-\lambda_3^2}+\frac{x_3^2}{\lambda_1^2-\lambda_2^2}+\frac{y_1^2}{\lambda_1^2-\lambda_4^2};$$
\medskip
$$G_2(M)=\frac{x_1^2}{\lambda_2^2-\lambda_3^2}+\frac{x_3^2}{\lambda_2^2-\lambda_1^2}+\frac{y_2^2}{\lambda_2^2-\lambda_4^2};$$
\medskip
$$G_3(M)=\frac{x_1^2}{\lambda_3^2-\lambda_2^2}+\frac{x_2^2}{\lambda_3^2-\lambda_1^2}+\frac{y_3^2}{\lambda_3^2-\lambda_4^2};$$
\medskip
$$G_4(M)=\frac{y_1^2}{\lambda_4^2-\lambda_1^2}+\frac{y_2^2}{\lambda_4^2-\lambda_2^2}+\frac{y_3^2}{\lambda_4^2-\lambda_3^2}.$$

\medskip

By direct computation we observe that $G_k$, $k=1,2,3,4$, are of the form \eqref{energy}. More precisely,
\begin{itemize}
\item $G_1=\mu_1\cdot H+ \mu_2\cdot I+m'\cdot C_1$, where
\begin{align*}
\mu_1&=\frac{2(\lambda_2 +\lambda_3 )(\lambda_2 +\lambda_4 )(\lambda_3 +\lambda_4 )}{(\lambda_1 +\lambda_2+\lambda_3 +\lambda_4 )(\lambda_1 -\lambda_ 2)(\lambda_1 -\lambda_3 )(\lambda_1 -\lambda_4 )};\\
\mu_2&=\frac{1}{(\lambda_1 +\lambda_2+\lambda_3 +\lambda_4 )(\lambda_1 -\lambda_ 2)(\lambda_1 -\lambda_3 )(\lambda_1 -\lambda_4 )};\\
m'&=-\frac{2(\lambda_2^2 +\lambda_3^2 +\lambda_4^2+\lambda_2 \lambda_3 +\lambda_2 \lambda_4 +\lambda_3 \lambda_4) }{(\lambda_1 +\lambda_2+\lambda_3 +\lambda_4 )(\lambda_1 -\lambda_ 2)(\lambda_1 -\lambda_3 )(\lambda_1 -\lambda_4 )};
\end{align*}

\item $G_2=\mu_1\cdot H+ \mu_2\cdot I+m'\cdot C_1$, where
\begin{align*}
\mu_1&=-\frac{2(\lambda_1 +\lambda_3 )(\lambda_1 +\lambda_4 )(\lambda_3 +\lambda_4 )}{(\lambda_1 +\lambda_2+\lambda_3 +\lambda_4 )(\lambda_1 -\lambda_ 2)(\lambda_2 -\lambda_3 )(\lambda_2 -\lambda_4 )};\\
\mu_2&=-\frac{1}{(\lambda_1 +\lambda_2+\lambda_3 +\lambda_4 )(\lambda_1 -\lambda_ 2)(\lambda_2 -\lambda_3 )(\lambda_2 -\lambda_4 )};\\
m'&=\frac{2(\lambda_1^2 +\lambda_3^2 +\lambda_4^2+\lambda_1 \lambda_3 +\lambda_1 \lambda_4 +\lambda_3 \lambda_4) }{(\lambda_1 +\lambda_2+\lambda_3 +\lambda_4 )(\lambda_1 -\lambda_ 2)(\lambda_2 -\lambda_3 )(\lambda_2 -\lambda_4 )};
\end{align*}

\item $G_3=\mu_1\cdot H+ \mu_2\cdot I+m'\cdot C_1$, where
\begin{align*}
\mu_1&=\frac{2(\lambda_1 +\lambda_2 )(\lambda_1 +\lambda_4 )(\lambda_2 +\lambda_4 )}{(\lambda_1 +\lambda_2+\lambda_3 +\lambda_4 )(\lambda_1 -\lambda_ 3)(\lambda_2 -\lambda_3 )(\lambda_3 -\lambda_4 )};\\
\mu_2&=\frac{1}{(\lambda_1 +\lambda_2+\lambda_3 +\lambda_4 )(\lambda_1 -\lambda_ 3)(\lambda_2 -\lambda_3 )(\lambda_2 -\lambda_4 )};\\
m'&=-\frac{2(\lambda_1^2 +\lambda_2^2 +\lambda_4^2+\lambda_1 \lambda_2 +\lambda_1 \lambda_4 +\lambda_2 \lambda_4) }{(\lambda_1 +\lambda_2+\lambda_3 +\lambda_4 )(\lambda_1 -\lambda_ 3)(\lambda_2 -\lambda_3 )(\lambda_3 -\lambda_4 )};
\end{align*}

\item $G_4=\mu_1\cdot H+ \mu_2\cdot I+m'\cdot C_1$, where
\begin{align*}
\mu_1&=-\frac{2(\lambda_1 +\lambda_2 )(\lambda_1 +\lambda_3 )(\lambda_2 +\lambda_3 )}{(\lambda_1 +\lambda_2+\lambda_3 +\lambda_4 )(\lambda_1 -\lambda_ 4)(\lambda_2 -\lambda_4 )(\lambda_3 -\lambda_4 )};\\
\mu_2&=-\frac{1}{(\lambda_1 +\lambda_2+\lambda_3 +\lambda_4 )(\lambda_1 -\lambda_ 4)(\lambda_2 -\lambda_4 )(\lambda_3 -\lambda_4 )};\\
m'&=\frac{2(\lambda_1^2 +\lambda_2^2 +\lambda_3^2+\lambda_1 \lambda_2 +\lambda_1 \lambda_3 +\lambda_2 \lambda_3) }{(\lambda_1 +\lambda_2+\lambda_3 +\lambda_4 )(\lambda_1 -\lambda_ 4)(\lambda_2 -\lambda_4 )(\lambda_3 -\lambda_4 )}.
\end{align*}
\end{itemize}

For the equilibrium $M_{a,b}^1$ in $\mathfrak{t}_1 \cap \operatorname{Orb}_{c_1; c_2}$ let us consider the Lyapunov function $F_{mn}^{\mu_1\mu_2}\in
C^{\infty}(\mathfrak{so}(4),\mathbb{R})$ 
$$F_{mn}^{\mu_1\mu_2}(M)=G_3(M)+m C_1(M)+n C_2(M),$$
where $m,n$ are real numbers.

Choosing $m,n$ such that $dF_{mn}^{\mu_1\mu_2}(M_{a,b}^1)=0$ and denoting
$W:=\ker dC_1(M_{a,b}^1)\cap \ker dC_2(M_{a,b}^1)$
we obtain the determinants associated with all upper-left submatrices of the
Hessian $d^2F_{mn}^{\mu_1\mu_2}(M_{a,b}^1)|_{W\times W}$
\begin{align*}
D_1&=2\cdot\frac{a^2(\lambda_1^2-\lambda_2^2)+b^2(\lambda_2^2-\lambda_3^2)}{(\lambda_1^2-\lambda_3^2)(\lambda_2^2-\lambda_3^2)(a^2-b^2)}>0;\\
D_2&=4a^2\cdot\frac{a^2(\lambda_1^2-\lambda_2^2)+b^2(\lambda_2^2-\lambda_3^2)}{(\lambda_1^2-\lambda_3^2)(\lambda_2^2-\lambda_3^2)^2(a^2-b^2)^2}>0;\\
D_3&=8a^4\cdot \frac{\lambda_1^2-\lambda_2^2}{(\lambda_1^2-\lambda_3^2)(\lambda_2^2-\lambda_3^2)^3(a^2-b^2)^2}>0;\\
D_4&=16a^4\cdot\frac{(\lambda_1^2-\lambda_2^2)(\lambda_2^2-\lambda_4^2)}{(\lambda_1^2-\lambda_3^2)(\lambda_3^2-\lambda_4^2)(\lambda_2^2-\lambda_3^2)^4(a^2-b^2)^2}>0.
\end{align*}
It follows that $d^2F_{mn}^{\mu_1\mu_2}(M_{a,b}^1)|_{W\times W}$ is positive definite and thus the equilibrium $M_{a,b}^1$ in $\mathfrak{t}_1 \cap \operatorname{Orb}_{c_1; c_2}$ is nonlinear stable. A similar computation proves that the equilibrium $M_{-a,-b}^1$ in $\mathfrak{t}_1 \cap \operatorname{Orb}_{c_1; c_2}$ is also nonlinear stable.

\medskip

Analogously, for equilibria $M_{a,b}^3$ and $M_{-a,-b}^3$ in $\mathfrak{t}_3 \cap \operatorname{Orb}_{c_1; c_2}$ a convenient Lyapunov function for applying Arnold's energy method is $F_{mn}^{\mu_1\mu_2}(M)=G_1(M)+m C_1(M)+n C_2(M)$, while for equilibria $M_{b,a}^3$ and $M_{-b,-a}^3$ in $\mathfrak{t}_3 \cap \operatorname{Orb}_{c_1; c_2}$ a convenient Lyapunov function for applying Arnold's energy method is $F_{mn}^{\mu_1\mu_2}(M)=G_4(M)+m C_1(M)+n C_2(M)$.

Thus, we have obtained the following result:
\begin{thm}
Using Lyapunov functions, we have the following stability results.
\begin{itemize}
\item[(i)] The equilibria $M_{a,b}^1,M_{-a,-b}^1$ are nonlinear stable with Lyapunov function
$$F_{mn}^{\mu_1\mu_2}(M)=G_3(M)+m C_1(M)+n C_2(M).$$
\item[(ii)] The equilibria $M_{a,b}^3,M_{-a,-b}^3$ are nonlinear stable with Lyapunov function
$$F_{mn}^{\mu_1\mu_2}(M)=G_1(M)+m C_1(M)+n C_2(M).$$
\item[(iii)] The equilibria $M_{b,a}^3,M_{-b,-a}^3$ are nonlinear stable with Lyapunov function
$$F_{mn}^{\mu_1\mu_2}(M)=G_4(M)+m C_1(M)+n C_2(M).$$
\end{itemize}
\end{thm}

\section{Spectral stability for equilibria with bifurcation behavior}

As we have seen in Theorem \ref{bifurcation}, for equilibria $M_{b,a}^1,M_{-b,-a}^1$ a bifurcation in the stability behavior does occur. In this section we aim to find the region of spectral stability for the equilibria $M_{b,a}^1$ and $M_{-b,-a}^1$. Notice that from Theorem \ref{values_a_b} we have $a\in (0,\infty),b\in [0,\infty)$ and $b<a$, which also implies that $0\leq \frac{b^2}{a^2}<1$.

Before we start the spectral stability analysis, we need a few considerations about the position on the real line of the roots $\alpha_1<\alpha_2$ of the function $\tilde{f}$ defined in Section 2. 
The discriminant of $\widetilde{f}$ is given by
$$\Delta_{\widetilde{f}}=16(\lambda_1^2-\lambda_2^2)(\lambda_1^2-\lambda_3^2)(\lambda_2^2-\lambda_4^2)(\lambda_3^2-\lambda_4^2)>0.$$
Since $\Delta_{\widetilde{f}}>0,-\frac{\widetilde{T}}{\widetilde{S}}>0$ and $\frac{\widetilde{U}}{\widetilde{S}}>0$, the quadratic equation associated to $\widetilde{f}$ has two distinct strictly positive real solutions. We notice that $\widetilde{f}(1)=(\lambda_1^2+\lambda_4^2- \lambda_2^2- \lambda_3^2)^2\geq 0$ and
$$-\frac{\widetilde{T}}{2\widetilde{S}}-1=-\frac{(\lambda_1^2-\lambda_2^2)(\lambda_1^2-\lambda_3^2)+(\lambda_2^2-\lambda_4^2)(\lambda_3^2-\lambda_4^2)}{(\lambda_1^2-\lambda_4^2)^2}<0,$$ which implies the following ordering
$$0<\alpha_1<\alpha_2\leq 1.$$
The above computation also shows that $\alpha_2\not= 1$ if and only if $\lambda_1^2+\lambda_4^2\not= \lambda_2^2+ \lambda_3^2$. 

The characteristic equation for the linearized system $L_{M_{b,a}^1}X_H|_{\operatorname{Orb}_{c_1;c_2}}$ is given by
\begin{equation}\label{delta}
ut^4+vt^2+w=0,
\end{equation}
where
$$u=(\lambda_1 + \lambda_2 )(\lambda_1 + \lambda_3 )(\lambda_1 + \lambda_4 )^4(\lambda_2 + \lambda_3 )^4(\lambda_2 + \lambda_ 4)(\lambda_3 + \lambda_ 4)>0;$$
$$w=(\lambda_1 - \lambda_2 )(\lambda_1 - \lambda_3 )(\lambda_2 - \lambda_4 )(\lambda_3 - \lambda_4 )[(\lambda_2 + \lambda_3 )^2a^2-(\lambda_1 + \lambda_4 )^2b^2]^2\geq 0$$
and 
$$v=-2(\lambda_1 + \lambda_4 )^2(\lambda_2 + \lambda_3 )^2a^2T\left(\frac{S}{T}-\frac{b^2}{a^2}\right).$$
We have made the notations
$$S=\left( \lambda_{{2}}+\lambda_{{3}} \right) ^{2} E_2;~~~
E_2=-(\lambda_2^2+\lambda_1\lambda_4)(\lambda_3^2+\lambda_1\lambda_4)+\lambda_2\lambda_3(\lambda_1 + \lambda_4 )^2$$
and respectively
$$T=\left( \lambda_{{1}}+\lambda_{{4}} \right) ^{2} E_1;~~~
E_1=(\lambda_1^2+\lambda_2\lambda_3)(\lambda_4^2+\lambda_2\lambda_3)-\lambda_1\lambda_4(\lambda_2 + \lambda_3 )^2.$$

Next, we will prove that $T>0$. Indeed, if $\lambda_4\leq0$, from the conditions $\lambda_i+\lambda_j>0$ we obtain $0<\lambda_3<\lambda_2<\lambda_1$ and consequently, $T>0$. If $\lambda_4>0$, then 
$\lambda_1^2+\lambda_2\lambda_3>\lambda_1(\lambda_2 + \lambda_3 )$ as it is equivalent with $(\lambda_1 - \lambda_2 )(\lambda_1 - \lambda_3 )>0$
and 
$\lambda_4^2+\lambda_2\lambda_3>\lambda_4(\lambda_2 + \lambda_3 )$ as it is equivalent with $(\lambda_4 - \lambda_2 )(\lambda_4 - \lambda_3 )>0$. By multiplication we obtain that $E_1>0$ and consequently, $T>0$.

In equation \eqref{delta} we denote $t^2=s$ and the discriminant of the corresponding quadratic equation is
$$\Delta=4(\lambda_1 + \lambda_ 4)^6(\lambda_2 + \lambda_3 )^6(\lambda_1\lambda_4-\lambda_2\lambda_3)^2a^4\tilde{f}\left(\frac{b^2}{a^2}\right).$$

Now, we study the position of $\frac{S}{T}$ with respect to the roots $\alpha_1,\alpha_2$. We notice that
$$\tilde{f}\left(\frac{S}{T}\right)=-\frac{(\lambda_1 ^2-\lambda_2 ^2)(\lambda_1 ^2-\lambda_3 ^2)(\lambda_2 ^2-\lambda_4 ^2)(\lambda_3 ^2-\lambda_4 ^2)
(\lambda_ 2+ \lambda_3 )^2(\lambda_1\lambda_4-\lambda_2\lambda_3)^2}{E_1^2(\lambda_1 + \lambda_4 )^2}\leq 0$$
and consequently $\alpha_1\leq \frac{S}{T}\leq \alpha_2$.

We discuss the spectral stability problem for $L_{M_{b,a}^1}X_H|_{\operatorname{Orb}_{c_1;c_2}}$ by analyzing the position of the number $\frac{b^2}{a^2}$ in the interval $[0,1)$. Taking into account the expression of $\Delta$, the following two cases: (I) $\lambda_1\lambda_4\not= \lambda_2\lambda_3$ and (II) $\lambda_1\lambda_4= \lambda_2\lambda_3$ arise naturally.\\\\
{\bf Case I ($\lambda_1\lambda_4\not= \lambda_2\lambda_3$)} \\
In this case $\tilde{f}\left(\frac{S}{T}\right)<0$ and consequently, $\alpha_1<\frac{S}{T}<\alpha_2$.\\\\
{\bf Subcase 1} $\frac{b^2}{a^2}\in [0,\alpha_1)$.\\
We have $\frac{b^2}{a^2}<\alpha_1<\frac{S}{T}$, which implies $\frac{S}{T}-\frac{b^2}{a^2}>0$ and consequently, $v<0$. Also $\Delta>0$ as we are outside the roots of $\tilde{f}$. It follows that the equation $us^2+vs+w=0$ has two distinct real positive roots $s_1<s_2$, where at least the root $s_2$ is strictly positive. Solving $t^2=s_2$, we obtain that equation \eqref{delta} has at least a solution with strictly positive real part, which implies {\bf spectral instability}.\\\\
{\bf Subcase 2} $\frac{b^2}{a^2}=\alpha_1$.\\
Then, $\frac{b^2}{a^2}=\alpha_1<\frac{S}{T}$, which implies $v<0$. Also $\Delta =0$. It follows that the equation $us^2+vs+w=0$ has a double real strictly positive root. Consequently, as before, we obtain {\bf spectral instability}.\\\\
{\bf Subcase 3} $\frac{b^2}{a^2}\in (\alpha_1,\alpha_2)$.\\
Then, $\Delta<0$. The characteristic equation \eqref{delta} has complex roots of the form $\pm A\pm \imath B$, with $A,B\in \mathbb{R}^*$. We obtain again {\bf spectral instability}.\\\\
{\bf Subcase 4} $\frac{b^2}{a^2}\in [\alpha_2,1)$.

If $\lambda_1^2+\lambda_4^2= \lambda_2^2+\lambda_3^2$, then $\alpha_2=1$ and consequently, the previous three subcases completely cover the interval $[0,1)$. 

Next, we analyze the situation when $\lambda_1^2+\lambda_4^2\not= \lambda_2^2+\lambda_3^2$. Then, $\alpha_2<1$ and we consider $\frac{b^2}{a^2}\in [\alpha_2,1)$.\\
a) $\frac{b^2}{a^2}=\alpha_2$. Then, $\Delta=0$. Also $\frac{S}{T}<\alpha_2=\frac{b^2}{a^2}$, which implies $\frac{S}{T}-\frac{b^2}{a^2}<0$ and consequently, $v>0$. The equation $us^2+vs+w=0$ has a double real strictly negative root. Then, equation \eqref{delta} has two double conjugate purely imaginary roots; we are in the case of {\bf spectral stability}.\\
b) $\frac{b^2}{a^2}\in (\alpha_2,1)$. Then, $\Delta>0$. Also $\frac{S}{T}<\alpha_2<\frac{b^2}{a^2}$, which implies $\frac{S}{T}-\frac{b^2}{a^2}<0$ and consequently, $v>0$. The equation $us^2+vs+w=0$ has two distinct negative real roots, at least one of them being strictly negative. Then, equation \eqref{delta} has roots of the form $\pm\imath A$ and $\pm\imath B$, $A\not= B,A>0, B\geq 0$. We are in the case of {\bf spectral stability}.\\\\
{\bf Case II ($\lambda_1\lambda_4= \lambda_2\lambda_3$)}\\
First of all, we notice that in this case $\lambda_1^2+\lambda_4^2\not= \lambda_2^2+\lambda_3^2$; otherwise, one would obtain $\lambda_1+\lambda_4= \lambda_2+\lambda_3$ and consequently, we obtain the set equality $\{\lambda_1, \lambda_4 \}=\{\lambda_2, \lambda_3 \}$. This, in turn, would contradict the initial hypothesis $\lambda_1> \lambda_2> \lambda_3> \lambda_4 $. Thus, for this case we have $\alpha_2\not= 1$. We have $\Delta=0$ and $\tilde{f}\left(\frac{S}{T}\right)=0$. By direct computation, $\alpha_1=\frac{\lambda_4^2(\lambda_2-\lambda_3)^2}{(\lambda_2\lambda_3-\lambda_4^2)^2}$ and $\alpha_2=\frac{S}{T}=\frac{\lambda_4^2(\lambda_2+\lambda_3)^2}{(\lambda_2\lambda_3-\lambda_4^2)^2}$.\\\\
{\bf Subcase 1} $\frac{b^2}{a^2}\in [0,\alpha_2)$. Then, $\frac{S}{T}=\alpha_2>\frac{b^2}{a^2}$, which implies $\frac{S}{T}-\frac{b^2}{a^2}>0$ and consequently, $v<0$. As in case I we obtain that equation \eqref{delta} has double roots of the form $\pm A$, with $A>0$. We are in the case of {\bf spectral instability}.\\\\
{\bf Subcase 2} $\frac{b^2}{a^2}=\alpha_2$. Then, $\frac{S}{T}=\frac{b^2}{a^2}$, which implies $v=0$. The equation $us^2+vs+w=0$ has $0$ as a double solution. Consequently, equation \eqref{delta} has $0$ as a root of multiplicity 4; we are in the case of {\bf spectral stability}.\\\\
{\bf Subcase 3} $\frac{b^2}{a^2}\in (\alpha_2,1)$. Then, $\alpha_2=\frac{S}{T}<\frac{b^2}{a^2}$, which implies $v>0$. The equation $us^2+vs+w=0$ has a double strictly negative real solution. Consequently, equation \eqref{delta} has double conjugate purely imaginary roots, which leads to {\bf spectral stability}.

\medskip

Identical results hold for the equilibrium $M_{-b,-a}^1$.
The above discussion can be summarized as follows:

\begin{thm}
\label{spectral}
The equilibria $M_{b,a}^1$ and $M_{-b,-a}^1$ have the following spectral stability behavior:\\\\
(i) if $\lambda_1^2+\lambda_4^2= \lambda_2^2+\lambda_3^2$, then $M_{b,a}^1$ and $M_{-b,-a}^1$ are spectrally unstable; if $\lambda_1^2+\lambda_4^2\not= \lambda_2^2+\lambda_3^2$, then $M_{b,a}^1$ and $M_{-b,-a}^1$ are spectrally unstable for $\frac{b^2}{a^2}\in [0,\alpha_2)$.\\\\
(ii) if $\lambda_1^2+\lambda_4^2\not= \lambda_2^2+\lambda_3^2$, then $M_{b,a}^1$ and $M_{-b,-a}^1$ are spectrally stable for $\frac{b^2}{a^2}\in [\alpha_2,1)$.

\end{thm}

We notice that this result is in agreement with Theorem \ref{bifurcation} and Theorem \ref{bifurcation2}.

\section{Nonlinear stability using energy methods for equilibria with bifurcation behavior}

For the equilibrium $M_{b,a}^1$, under the condition $\lambda_1^2+\lambda_4^2\not= \lambda_2^2+\lambda_3^2$, we have found spectral stability for $\frac{b^2}{a^2}\in [\alpha_2,1)$. For the equilibrium $M_{b,a}^1$ in $\mathfrak{t}_1\cap \operatorname{Orb}_{c_1;c_2}$ 
it is natural to choose a Lyapunov function of the form $F_{mn}^{\mu_1\,u_2}=G_2+mC_1+nC_2$. Applying Arnold's method we find indefiniteness. We will introduce an additional degree of freedom and thus we search for a Lyapunov function of the following form:
$$F_{mn}^p=H+\frac{1}{2(\lambda_1 +\lambda_3 )(\lambda_1 +\lambda_4 )(\lambda_3 +\lambda_4 )}\cdot p\cdot I+m\cdot C_1+n\cdot C_2.$$
The coefficient 
$\displaystyle\frac{1}{2(\lambda_1 +\lambda_3 )(\lambda_1 +\lambda_4 )(\lambda_3 +\lambda_4 )}$ is equal with $\frac{\mu_2}{\mu_1},$
where $\mu_1,\mu_2$ are the coefficients that appear in the expression of $G_2$.\\
Choosing $m,n$ such that $dF_{mn}^p(M_{b,a}^1)=0$ and denoting
$W:=\ker dC_1(M_{b,a}^1)\cap \ker dC_2(M_{b,a}^1)$,
we obtain the determinants associated with all upper-left submatrices of the
Hessian $d^2F_{mn}^p(M_{b,a}^1)|_{W\times W}$:
\begin{itemize}
\item $D_1=\displaystyle \frac{(S_1a^2-T_1b^2)p+S_1'a^2+T_1'b^2}{
(\lambda_1 + \lambda_3 )(\lambda_1 + \lambda_4 )(\lambda_2 + \lambda_3 )(\lambda_ 3+ \lambda_4 )(a^2-b^2)}.$\\
We denote
$D_1':=(S_1a^2-T_1b^2)p+S_1'a^2+T_1'b^2$, where:
\begin{align*}
S_1&=(\lambda_2 + \lambda_3 )(\lambda_3 ^2-\lambda_4 ^2)>0;
T_1=(\lambda_2 + \lambda_3 )(\lambda_1 ^2-\lambda_2 ^2)>0;\\
S_1'&=-(\lambda_2 + \lambda_3 )(\lambda_3 ^2-\lambda_4 ^2);
T_1'=(\lambda_1 - \lambda_2 )(\lambda_1 + \lambda_4 )(\lambda_3 + \lambda_4 ).\end{align*}
\item $D_2=\displaystyle\frac{D_1'\cdot [(S_2a^2-T_2b^2)p+S_2'a^2+T_2'b^2]}{(\lambda_1 + \lambda_2 )
(\lambda_1 + \lambda_3 )^2(\lambda_1 + \lambda_4 )^2(\lambda_2 + \lambda_3 )^2(\lambda_ 3+ \lambda_4 )^2(a^2-b^2)^2}.$\\
We denote
$D_2':=(S_2a^2-T_2b^2)p+S_2'a^2+T_2'b^2$,
where:
\begin{align*}
S_2&=(\lambda_1 + \lambda_2 )(\lambda_2 + \lambda_3 )(\lambda_2 ^2-\lambda_4 ^2)>0;
T_2=(\lambda_1 + \lambda_2 )(\lambda_2 + \lambda_3 )(\lambda_1 ^2-\lambda_3 ^2)>0;\\
S_2'&=-(\lambda_1 + \lambda_3 )(\lambda_2 + \lambda_3 )(\lambda_3 + \lambda_4 )(\lambda_2 - \lambda_4 );
T_2'=(\lambda_1 + \lambda_3 )(\lambda_1 + \lambda_4 )(\lambda_ 3+ \lambda_4 )(\lambda_1 - \lambda_3 ).\end{align*}
\item $D_3=-\displaystyle\frac{D_2'(S_3p^2+T_3p+U_3)(\lambda_ 1- \lambda_2 )(\lambda_ 3- \lambda_4 )}{(\lambda_1 + \lambda_2 )
(\lambda_1 + \lambda_3 )^3(\lambda_1 + \lambda_4 )^3(\lambda_2 + \lambda_3 )^3(\lambda_2 + \lambda_4 )(\lambda_ 3+ \lambda_4 )^2(a^2-b^2)^2}$.\\
We denote $g(p):=S_3p^2+T_3p+U_3$, where:
\begin{align*}
S_3&=(\lambda_1 + \lambda_2 )(\lambda_2 + \lambda_ 4)(\lambda_2 + \lambda_3 )^2(a^2-b^2)>0;\\ 
T_3&=-2(\lambda_ 2+\lambda_3 )\{(\lambda_2 +\lambda_ 3)[(\lambda_1 +\lambda_2 )(\lambda_2 +\lambda_4 )+(\lambda_1 +\lambda_3 )(\lambda_3 +\lambda_4 )]a^2-\\
 &-(\lambda_ 1+\lambda_ 4)[(\lambda_ 1+\lambda_ 2)(\lambda_1 +\lambda_3 )+(\lambda_2 +\lambda_4 )(\lambda_3 +\lambda_4 )]b^2\};\\
U_3&=(\lambda_1 +\lambda_3 )(\lambda_3 +\lambda_4 )[(\lambda_2 + \lambda_3 )^2a^2-(\lambda_1 + \lambda_4 )^2b^2].\end{align*}
\item $D_4=\displaystyle\frac{(\lambda_1 - \lambda_2 )(\lambda_1 - \lambda_3 )(\lambda_2 - \lambda_4 )(\lambda_3 - \lambda_4 )(S_3p^2+T_3p+U_3)^2}{(\lambda_1 + \lambda_2 )
(\lambda_1 + \lambda_3 )^3(\lambda_1 + \lambda_4 )^4(\lambda_2 + \lambda_3 )^4(\lambda_2 + \lambda_4 )(\lambda_ 3+ \lambda_4 )^3(a^2-b^2)^2}.$
\end{itemize}

We are looking for $p\in \mathbb{R}$ such that Hess$F_{mn}^p(M_{b,a}^1)|_{W\times W}$ is definite. This is equivalent with $D_1D_3>0,D_2>0,D_4>0$. We notice that $D_1D_3$ has the sign of $-D_1'D_2'g(p)$. Since $D_2$ has the sign of $D_1'D_2'$, for definiteness we need to have $g(p)<0$. Thus, definiteness is equivalent with the existence of a real number $p$ such that
\begin{equation}
\label{definite}
D_1'(p)D_2'(p)>0,~~~g(p)<0.
\end{equation}

We distinguish two cases: (I) $\lambda_1^2+\lambda_4^2< \lambda_2^2+\lambda_3^2$ and (II) $\lambda_1^2+\lambda_4^2> \lambda_2^2+\lambda_3^2$.

\medskip

{\bf Case I $\lambda_1^2+\lambda_4^2< \lambda_2^2+\lambda_3^2$} \\
We notice that 
$\frac{S_1}{T_1}=\frac{\lambda_3^2-\lambda_4^2}{\lambda_1^2-\lambda_2^2}>1>\frac{b^2}{a^2}$ and $\frac{S_2}{T_2}=\frac{\lambda_2^2-\lambda_4^2}{\lambda_1^2-\lambda_3^2}>1>\frac{b^2}{a^2}$, thus $\frac{S_1}{T_1}-\frac{b^2}{a^2}>0$ and $\frac{S_2}{T_2}-\frac{b^2}{a^2}>0$. It follows that the linear functions $D_1'(p),D_2'(p)$, as functions of variable $p$, are strictly increasing and the associated equations $D_1'(p)=0$ and $D_2'(p)=0$ have each a unique root, which we denote by $p_1$, respectively $p_2$. We notice that
$$p_1-p_2=\frac{(\lambda_2 - \lambda_3 )(\lambda_2 - \lambda_ 4)(\lambda_3^2 - \lambda_4^2 )(\lambda_1+\lambda_2+\lambda_3+\lambda_4)(a^2-b^2)}{a^2(\lambda_1^2-\lambda_2^2)(\lambda_1^2-\lambda_3^2)(\lambda_1 + \lambda_2 )\left(\frac{S_1}{T_1}-\frac{b^2}{a^2}\right)\left(\frac{S_2}{T_2}-\frac{b^2}{a^2}\right)}>0,$$
which implies $p_1>p_2$.

By direct computation we obtain that the discriminant of the quadratic equation $g(p)=0$ is
$$\Delta_3=T_3^2-4S_3U_3=(\lambda_1+\lambda_2+\lambda_3+\lambda_4)^2(\lambda_2 + \lambda_3 )^2a^4\tilde{f}\left(\frac{b^2}{a^2}\right).$$
For $\frac{b^2}{a^2}=\alpha_2$ we have that $\Delta_3=0$ and consequently \eqref{definite} does not have a solution.\\
For $\frac{b^2}{a^2}>\alpha_2$ we have $\Delta_3>0$ and the equation $g(p)=0$ has two distinct real roots, which we denote by $p_3<p_4$.

In order to solve the system of inequalities \eqref{definite} we need to study the position of $p_1,p_2$ with respect to $p_3,p_4$. For this, we need the following computations:
$$g(p_1)=\frac{b^2(\lambda_1 - \lambda_2 )(\lambda_2 - \lambda_4 )^2
(\lambda_ 3^2-\lambda_ 4^2)(\lambda_1+\lambda_2+\lambda_3+\lambda_4)^2(a^2-b^2)}{a^2 (\lambda_1 + \lambda_2 )(\lambda_1^2 - \lambda_2^2 )\left(\frac{S_1}{T_1}-\frac{b^2}{a^2}\right)^2}>0$$
and
$$g(p_2)=\frac{b^2(\lambda_2 - \lambda_4 )
(\lambda_ 3^2-\lambda_ 4^2)^2(\lambda_1+\lambda_2+\lambda_3+\lambda_4)^2(a^2-b^2)}{a^2(\lambda_1 + \lambda_2 )(\lambda_1^2 - \lambda_3^2 ) \left(\frac{S_2}{T_2}-\frac{b^2}{a^2}\right)^2}>0.$$
Consequently, $p_1$ and $p_2$ are outside of the interval $(p_3,p_4)$.

Next, in order to compare $p_1$ with the minimum point $-\frac{T_3}{2S_3}$ of the parabola $g(p)$, we prove:
\begin{lem}
If $\lambda_1^2+\lambda_4^2< \lambda_2^2+\lambda_3^2$, then the following inequality holds:
$$p_1-\left(-\frac{T_3}{2S_3}\right)<0.$$
\end{lem}
\noindent {\bf Proof.} By direct computation we obtain
$$p_1-\left(-\frac{T_3}{2S_3}\right)=\frac{1}{2}\cdot \frac{(\lambda_1+\lambda_2+\lambda_3+\lambda_4)b^4\left[S_4\left(\frac{a^2}{b^2}\right)^2+T_4\frac{a^2}{b^2}+U_4\right]}{(\lambda_1 + \lambda_2 )(\lambda_2 + \lambda_3 )(\lambda_2 + \lambda_4 )(a^2-b^2)a^2(\lambda_1^2 - \lambda_2^2)\left(\frac{S_1}{T_1}-\frac{b^2}{a^2}\right)},$$
where:
\begin{align*}
S_4&=(\lambda_2 ^2-\lambda_3 ^2)(\lambda_3 ^2-\lambda_4 ^2)>0;\\
T_4&=2\,{\lambda_{{1}}}^{2}{\lambda_{{3}}}^{2}-3\,{\lambda_{{4}}}^{2}{
\lambda_{{1}}}^{2}+{\lambda_{{1}}}^{2}{\lambda_{{2}}}^{2}+4\,{\lambda_
{{2}}}^{2}{\lambda_{{4}}}^{2}-{\lambda_{{2}}}^{4}+{\lambda_{{4}}}^{2}{
\lambda_{{3}}}^{2}-3\,{\lambda_{{2}}}^{2}{\lambda_{{3}}}^{2}-{\lambda_
{{4}}}^{4};\\
U_4&=-(\lambda_1^2-\lambda_2^2)(\lambda_1^2-\lambda_4^2).
\end{align*}

We denote
$h_1(t):=\tilde{U}t^2+\tilde{T}t+\tilde{S}$ and
$h_2(t):=S_4t^2+T_4t+U_4$. Obviously, the quadratic equation $h_1(t)=0$ has the roots $\frac{1}{\alpha_2}<\frac{1}{\alpha_1}$. 
We notice that the quadratic function $h_3(t):=h_1(t)+h_2(t)$ has the distinct real roots:
$$\theta_1=1;~~~\theta_2=\frac{\lambda_1^2-\lambda_4^2}{\lambda_2^2-\lambda_3^2}>1.$$
It follows that
$$h_3(t)=(S_4+\tilde{U})(t-1)\left(t-\frac{\lambda_1^2-\lambda_4^2}{\lambda_2^2-\lambda_3^2}\right).$$
It is easy to see that $S_4+\tilde{U}=(\lambda_2^2-\lambda_3^2)(\lambda_2^2-\lambda_4^2)>0$. 
Also,
$$h_1\left(\frac{\lambda_1^2-\lambda_4^2}{\lambda_2^2-\lambda_3^2}\right)=-\frac{4(\lambda_1^2-\lambda_2^2)(\lambda_1^2-\lambda_4^2)(\lambda_3^2-\lambda_4^2)}{\lambda_2^2-\lambda_3^2}<0,$$
which implies that $\frac{1}{\alpha_2}< \frac{\lambda_1^2-\lambda_4^2}{\lambda_2^2-\lambda_3^2} <\frac{1}{\alpha_1}$.

Consequently, we have that $\left(1,\frac{1}{\alpha_2}\right)\subset \left(1,\frac{\lambda_1^2-\lambda_4^2}{\lambda_2^2-\lambda_3^2}\right)$, which implies that for $t\in \left(1,\frac{1}{\alpha_2}\right)$ we have that $h_3(t)<0$, which leads to $h_2(t)<-h_1(t)<0$. Since $\frac{a^2}{b^2}\in \left(1,\frac{1}{\alpha_2}\right)$ one obtains $h_2\left(\frac{a^2}{b^2}\right)<0$. This implies $p_1-\left(-\frac{T_3}{2S_3}\right)<0$. 
\rule{0.5em}{0.5em}

The above Lemma implies the following ordering
$$p_2<p_1<p_3<p_4.$$

The system of inequalities \eqref{definite} has as solutions any $p\in (p_3,p_4)$; for such $p$ we have $D_1'(p)>0$, $D_2'(p)>0$ and $g(p)<0$.
In conclusion, there exists a real number $p$ such that Hess$F_{mn}^p(M_{b,a}^1)|_{W\times W}$ is {\bf positive definite}.

\medskip

{\bf Case II $\lambda_1^2+\lambda_4^2> \lambda_2^2+\lambda_3^2$}\\\\
First of all, as we did in Case I, we need to establish the signs of the expressions $\frac{S_1}{T_1}-\frac{b^2}{a^2}$ and $\frac{S_2}{T_2}-\frac{b^2}{a^2}$. By direct computations we obtain
$$\tilde{f}\left(\frac{S_2}{T_2}\right)=-\frac{(\lambda_3 ^2-\lambda_4 ^2)(\lambda_1^2+\lambda_4^2-\lambda_2^2-\lambda_3^2)\cdot E}{(\lambda_1 ^2-\lambda_ 3^2)^2},$$
where
$$E=-3\,{\lambda_{{2}}}^{2}{\lambda_{{3}}}^{2}-{\lambda_{{2}}}^{2}{\lambda
_{{4}}}^{2}+4\,{\lambda_{{1}}}^{2}{\lambda_{{2}}}^{2}+{\lambda_{{3}}}^
{4}-3\,{\lambda_{{4}}}^{2}{\lambda_{{1}}}^{2}-{\lambda_{{1}}}^{2}{
\lambda_{{3}}}^{2}+2\,{\lambda_{{4}}}^{2}{\lambda_{{3}}}^{2}+{\lambda_
{{4}}}^{4}.$$
In order to decide the sign of the expression $E$, we make the following notations:
$$A_2:=\lambda_1 ^2-\lambda_2 ^2,~~~A_3:=\lambda_1 ^2-\lambda_ 3^2,~~~A_4:=\lambda_1 ^2-\lambda_4 ^2.$$
From the ordering $\lambda_1>\lambda_2>\lambda_3>\lambda_4$ we obtain that
$A_4>A_3>A_2>0$. We have
$$E=2A_3(A_4-A_2)+A_4(A_4-A_2)+A_3(A_3-A_2)>0.$$

Consequently, $\tilde{f}\left(\frac{S_2}{T_2}\right)<0$, thus $\frac{S_2}{T_2}<\alpha_2$. Since $\frac{b^2}{a^2}\in [\alpha_2,1)$, it follows that $\frac{S_2}{T_2}-\frac{b^2}{a^2}<0$.

We also have 
$$\frac{S_1}{T_1}-\frac{S_2}{T_2}=-\frac{(\lambda_2 ^2-\lambda_3 ^2)(\lambda_1^2+\lambda_4^2-\lambda_2^2-\lambda_3^2)}{(\lambda_1 ^2-\lambda_2 ^2)(\lambda_ 1^2-\lambda_3 ^2)}<0,$$
which implies
$\frac{S_1}{T_1} < \frac{S_2}{T_2} <\frac{b^2}{a^2}$ and consequently, $\frac{S_1}{T_1}-\frac{b^2}{a^2}<0$.

It follows that the linear functions $D_1'(p),D_2'(p)$, as functions of variable $p$, are strictly decreasing and the associated equations $D_1'(p)=0$ and $D_2'(p)=0$ have each a unique root, which we denote by $p_1$, respectively $p_2$. We notice that
$$p_1-p_2=\frac{(\lambda_2 - \lambda_3 )(\lambda_2 - \lambda_ 4)(\lambda_3^2 - \lambda_4^2 )(\lambda_1+\lambda_2+\lambda_3+\lambda_4)(a^2-b^2)}{a^2(\lambda_1^2-\lambda_2^2)(\lambda_1^2-\lambda_3^2)(\lambda_1 + \lambda_2 )\left(\frac{S_1}{T_1}-\frac{b^2}{a^2}\right)\left(\frac{S_2}{T_2}-\frac{b^2}{a^2}\right)}>0,$$
which implies $p_1>p_2$.

By direct computation we obtain that the discriminant of the quadratic equation $g(p)=0$ is
$$\Delta_3=T_3^2-4S_3U_3=(\lambda_1+\lambda_2+\lambda_3+\lambda_4)^2(\lambda_2 + \lambda_3 )^2a^4\tilde{f}\left(\frac{b^2}{a^2}\right).$$
For $\frac{b^2}{a^2}=\alpha_2$ we have that $\Delta_3=0$ and consequently \eqref{definite} does not have a solution.\\
For $\frac{b^2}{a^2}>\alpha_2$ we have $\Delta_3>0$ and the equation $g(p)=0$ has two distinct real roots, which we denote by $p_3<p_4$.

In order to solve the system of inequalities \eqref{definite} we need to study the position of $p_1,p_2$ with respect to $p_3,p_4$. For this, we need the following computations:
$$g(p_1)=\frac{b^2(\lambda_1 - \lambda_2 )(\lambda_2 - \lambda_4 )^2
(\lambda_ 3^2-\lambda_ 4^2)(\lambda_1+\lambda_2+\lambda_3+\lambda_4)^2(a^2-b^2)}{a^2 (\lambda_1 + \lambda_2 )(\lambda_1^2 - \lambda_2^2 )\left(\frac{S_1}{T_1}-\frac{b^2}{a^2}\right)^2}>0$$
and
$$g(p_2)=\frac{b^2(\lambda_2 - \lambda_4 )
(\lambda_ 3^2-\lambda_ 4^2)^2(\lambda_1+\lambda_2+\lambda_3+\lambda_4)^2(a^2-b^2)}{a^2(\lambda_1 + \lambda_2 )(\lambda_1^2 - \lambda_3^2 ) \left(\frac{S_2}{T_2}-\frac{b^2}{a^2}\right)^2}>0.$$
Consequently, $p_1$ and $p_2$ are outside of the interval $(p_3,p_4)$.

Next, in order to compare $p_2$ with the minimum point $-\frac{T_3}{2S_3}$ of the parabola $g(p)$, we prove
\begin{lem}
If $\lambda_1^2+\lambda_4^2> \lambda_2^2+\lambda_3^2$, then the following inequality holds:
$$p_2-\left(-\frac{T_3}{2S_3}\right)>0.$$
\end{lem}
\noindent {\bf Proof.} By direct computation we obtain
$$p_2-\left(-\frac{T_3}{2S_3}\right)=-\frac{1}{2}\cdot \frac{(\lambda_1+\lambda_2+\lambda_3+\lambda_4)b^4\left[S_5\left(\frac{a^2}{b^2}\right)^2+T_5\frac{a^2}{b^2}+U_5\right]}{(\lambda_1 + \lambda_2 )(\lambda_2 + \lambda_3 )(\lambda_2 + \lambda_4 )(a^2-b^2)a^2(\lambda_1^2 - \lambda_3^2)\left(\frac{S_2}{T_2}-\frac{b^2}{a^2}\right)},$$
where
\begin{align*}
S_5&=(\lambda_2^2-\lambda_3^2)(\lambda_2^2-\lambda_4^2);\\
T_5&=3\,{\lambda_{{4}}}^{2}{\lambda_{{1}}}^{2}-{\lambda_{{1}}}^{2}{\lambda_
{{3}}}^{2}-2\,{\lambda_{{1}}}^{2}{\lambda_{{2}}}^{2}+{\lambda_{{3}}}^{
4}-4\,{\lambda_{{4}}}^{2}{\lambda_{{3}}}^{2}+{\lambda_{{4}}}^{4}+3\,{
\lambda_{{2}}}^{2}{\lambda_{{3}}}^{2}-{\lambda_{{2}}}^{2}{\lambda_{{4}
}}^{2};\\
U_5&=(\lambda_1^2-\lambda_3^2)(\lambda_1^2-\lambda_4^2).
\end{align*}
We denote by
$h_4(t):=S_5t^2+T_5t+U_5$. We notice that the quadratic function $h_5(t):=h_4(t)-h_1(t)$ has the distinct real roots:
$$\theta_3=-\frac{\lambda_1^2-\lambda_4^2}{\lambda_2^2-\lambda_3^2}<0;~~~\theta_4=1.$$
It follows that
$$h_5(t)=(S_5-\tilde{U})\left(t+\frac{\lambda_1^2-\lambda_4^2}{\lambda_2^2-\lambda_3^2}\right)(t-1).$$
It is easy to see that $S_5-\tilde{U}=(\lambda_2^2-\lambda_3^2)(\lambda_3^2-\lambda_4^2)>0$. Consequently, for $t>1$ we have that $h_5(t)>0$, which is equivalent to $h_4(t)>h_1(t)$ for all $t>1$. Since $\frac{a^2}{b^2}\in \left(1,\frac{1}{\alpha_2}\right)$ one obtains $h_4\left(\frac{a^2}{b^2}\right)>
h_1\left(\frac{a^2}{b^2}\right)>0$. This implies $p_2-\left(-\frac{T_3}{2S_3}\right)>0$. \rule{0.5em}{0.5em}

The above Lemma implies the following ordering
$$p_3<p_4<p_2<p_1.$$

The system of inequalities \eqref{definite} has as solutions any $p\in (p_3,p_4)$. More precisely, for $p\in (p_3,p_4)$ we have $D_1'(p)>0$, $D_2'(p)>0$ and $g(p)<0$.
In conclusion, there exists a real number $p$ such that Hess$F_{mn}^p(M_{b,a}^1)|_{W\times W}$ is {\bf positive definite}.\\

The nonlinear stability for equilibria $M_{b,a}^1$ and $M_{-b,-a}^1$ can be summarized as follows:
\begin{thm}

We have the following stability behavior on regular orbits $\operatorname{Orb}_{c_1;c_2}$:
\begin{itemize}
\item[(i)] for $\frac{b^2}{a^2}\in (\alpha_2,1)$ the equilibria $M_{b,a}^1$ and $M_{-b,-a}^1$ are nonlinear stable with Lyapunov function
$$F_{mn}^p=H+\frac{p}{2(\lambda_1 +\lambda_3 )(\lambda_1 +\lambda_4 )(\lambda_3 +\lambda_4 )}I+mC_1+nC_2,$$
where $p\in (p_3,p_4)$;
\item[(ii)] for $\frac{b^2}{a^2}=\alpha_2$ the stability problem for the equilibria $M_{b,a}^1$ and $M_{-b,-a}^1$ cannot be decided using the set of constants of motion $\{H,I,C_1,C_2\}$; 
\item[(iii)] for $\frac{b^2}{a^2}\in [0,\alpha_2)$ the equilibria $M_{b,a}^1$ and $M_{-b,-a}^1$ are unstable.
\end{itemize}
\end{thm}

\medskip

\noindent {\bf Acknowledgments.} This work was supported by a grant of the Romanian National Authority for Scientific Research, CNCS – UEFISCDI, project number PN-II-RU-TE-2011-3-0006.

\end{document}